\documentclass[11pt,reqno,twoside]{article}
\usepackage{amsmath,amssymb,theorem,color,epsfig,epic,subfigure,amsfonts,hhline,graphicx}
\usepackage[normalem]{ulem}
\usepackage[colorinlistoftodos]{todonotes}

\theorembodyfont{\itshape}
\newtheorem{thm}{Theorem}[section]
\newtheorem{lem}[thm]{Lemma}
\newtheorem{prop}[thm]{Proposition}
{\theorembodyfont{\upshape}
\newtheorem{define}[thm]{Definition}
\newtheorem{rem}[thm]{Remark}

}

\newcommand{\Proof}[1][]{\noindent{\itshape Proof#1. }}
\newcommand{\EndProof}{\hfill$\Box$\bigskip}

\makeatletter
\let\@fnsymbol\@arabic
    
    \@addtoreset{equation}{section}
\makeatother

\newcommand\unit{\hbox{\rm 1\kern-2.8truept l}}

\newcommand\Lform{{\mathcal{L}}\kern-7.56pt\raise1.0pt\hbox{$-$}}

\begin{document}
\title{Quasi-invariant states with uniformly bounded cocycles}
\author{ Ameur Dhahri\footnote{Dipartimento di Matematica, Politecnico di Milano, Piazza Leonardo di Vinci 32, I-20133 Milano, Italy. E-mail: ameur.dhahri@polimi.it}, \'Eric Ricard \footnote{UNICAEN, CNRS, LMNO, 14000 Caen, France. E-mail: 
eric.ricard@unicaen.fr}}
\date{ }
  \maketitle

\begin{abstract}
We investigate the notion of quasi-invariant states introduced in \cite{AcDh1} from an analytic viewpoint. We give the structures of quasi-invariant states with uniformly bounded cocycles. As a consequence, we can apply a Theorem of Kov\'acs and Sz\"ucs to get a conditional expectation on fixed points and another of St{$\o$}rmer to get an invariant semifinite trace under extra assumptions.

\end{abstract}

\section{Introduction}
Let $(S,\mathcal B)$ be a measurable space and $G$ a group of $*$-automorphisms acting on $S$ by measurable transformation
$$s\in S\mapsto gs\in S,\quad g\in G.$$
A measure $\mu$ on $(S,\mathcal B)$ is called $G$-quasi-invariant if
the measures $\mu\circ g$ and $\mu$ are absolutely continuous one with respect to the other.

Quasi-invariant measures are a much wider class than invariant, and that is why they constitute the natural environment for the theory
of dynamical systems and ergodic theory. Therefore to develop a real quantum analogue of the classical theory of
dynamical systems, requires developing a theory of quasi-invariant states on general operator algebras. Recently, in \cite{AcDh1} the authors introduced the theory of quasi-invariant states under the action of a group of $*$-automorphisms on a von Neumann algebra (or a $*$-algebra). The relationship between the group $G$ and modular automorphism group, invariant subalgebras, ergodicity, modular theory, and abelian subalgebras is studied in \cite{DKY1}. Moreover, a martingale convergence theorem for strongly quasi-invariant states is given in \cite{DKY2}. Finally, in \cite{NT},  quasi-invariant actions are used in the setting of discrete quantum groups. In particular, the main technical results in \cite{NT} is the computation of the Radon-Nikodym cocycle for the action of the dual of $SUq(2)$ on its Martin boundary with respect to the canonical harmonic state defined by a quantum random walk.

In this note, we wish to clarify the notion of $G$-quasi-invariant state under a natural boundedness assumption. 
The paper is structured as follows. First, we give the structure of quasi-invariant states with uniformly bounded cocycle. Moreover, for $G$-quasi-invariant states, we show that one can always spatially implement the action of $G$. We further deal with a Theorem of Kov\'acs and Sz\"ucs. Then, we extend the main result of St\o rmer in \cite{St} to strongly quasi-invariant states with uniformly bounded cocycles.

 \section{Quasi-invariant states}\label{sec:general_theories}
We briefly review some definitions and properties of quasi-invariant states, respectively, strongly quasi-invariant states with respect to a group of $*$-automorphisms on a von Neumann algebra (cf \cite{AcDh1}). 

Let $\mathcal A$ be a von Neumann algebra and $G$ be a group of normal $*$-automorphisms of $\mathcal A$ equipped with its usual strong-topology. 
\begin{define} \label{def:quasi_invariance}
A faithful normal state $\varphi$ on $\mathcal A$ is said to be $G$-quasi-invariant if for all $g\in G$ there exists $x_g\in \mathcal A$ such that
\begin{equation}\label{eq:quasi_invariance}
\varphi(g(a))=\varphi(x_ga),\quad a\in \mathcal A.
\end{equation}
We say that $\varphi$ is $G$-strongly quasi-invariant if it is $G$-quasi-invariant and the Radon-Nikodym derivatives $x_g$ are self-adjoint: $x_g=x_g^*$, $g\in G$.
\end{define}

It is proved in \cite{AcDh1} the following properties:
\begin{enumerate}
 \item[(i)] When $\varphi$ is $G$-quasi-invariant, the map $g\mapsto x_g$ is a normalized multiplicative left $G$-1-cocycle satisfying
\begin{equation}\label{eq:cocycle}
x_e=\unit,\quad x_{g_2g_1}=x_{g_1}g_1^{-1}(x_{g_2}), \quad g_1,g_2\in G,
\end{equation}
and $x_g$ is invertible with its inverse
\begin{equation}\label{xg-1}
x_g^{-1}=g^{-1}(x_{g^{-1}}).
\end{equation}
Furthermore, for all $a\in \mathcal A$ and $g\in G$ it holds
\begin{equation}\label{xg*}
\varphi(x_ga)=\varphi(ax_g^*).
\end{equation}
\item[(ii)] If $\varphi$ is $G$-strongly quasi-invariant, it holds that $x_g$ is positive invertible. It also holds that $x_g$ commutes with $x_h$ for all $g,h\in G$. Hence the von Neumann algebra $\mathcal C$ generated by $\{x_g:g\in G\}$ is commutative. It can be also easily shown that $\mathcal C$ is a subalgebra of $\mathrm{Centr}(\varphi)$, the centralizer of $\varphi$, which is defined  by
\begin{equation}\label{eq:centralizer}
\mathrm{Centr}(\varphi):=\{x\in \mathcal A:\varphi(xy)=\varphi(yx)\text{ for all }y\in \mathcal A\}.
\end{equation}
\end{enumerate}

Now we introduce the following lemma  which plays a crucial rule later.
\begin{lem}[\cite{SZ}, 5.20. Lemma] \label{SZ} \rm{
Let $\varphi$ be a positive form on a $C^*$-algebra $\mathcal A$ and $a\in\mathcal A$. If the linear form $L_a\varphi$ defined on $\mathcal A$ by $L_a\varphi(x)=\varphi(ax)$ is positive, then 
$$L_a\varphi(x)\leq ||a||\varphi(x), $$
for all positive element $x\in\mathcal A$.
}
\end{lem}

\begin{rem}\label{abscont}
If $\mathcal A$ is a von Neumann algebra and $\varphi$, $\psi$ are two states on $\mathcal A$, we say that $\varphi$ is $\psi$--absolutely continuous if $\psi(a^*a)=0$ implies to $\varphi(a^*a)=0$ for $a\in \mathcal A$. We say that $\varphi$ and $\psi$ are equivalent ($\varphi \sim\psi$) if they are absolutely continuous one with respect to the other. A state $\varphi$ is  $G$-invariant (resp. $G$--quasi-inavriant) if $\varphi\circ g=\varphi$ (resp. $\varphi\circ g\sim \varphi$) for all $g\in G$. Indeed, if $\varphi$ is as in Definition \ref{def:quasi_invariance}, then by Lemma \ref{SZ} one has
$$\varphi\circ g (a^*a)=\varphi(x_ga^*a)\leq ||x_g||\varphi(a^*a),\quad \varphi(a^*a)\leq ||x_g^{-1}||\varphi\circ g(a^*a), \quad \forall g\in G,\;\forall a\in \mathcal A, $$
which proves that $\varphi\circ g\sim \varphi$ for all $g\in G$.

\end{rem}

\section{Quasi-invariant states with uniformly bounded cocycles}

We will make the following assumption.

{\bf Assumption}. Assume that there exists $\lambda>0$ such that 
\begin{equation}\tag{A}\label{A}||x_g||\leq \lambda,\quad\forall g\in G. \end{equation}
Since $x_e=\unit$, it is clear that $\lambda\geq 1$. By \eqref{xg-1}, we also have
$\| x_g ^{-1}\|\leq \lambda$.


\begin{lem}\label{lemm}{\rm
Under assumption \eqref{A}, for any $g\in G$ and any positive operator $a\in\mathcal A$}
$$\frac{1}{\lambda}\varphi(a)\leq \varphi(x_g a)\leq \lambda \varphi(a) \quad \textrm{and} \quad \frac{1}{\lambda}\varphi(a)\leq \varphi(a (x_g^{-1})^* )\leq \lambda \varphi(a).$$
\end{lem}
\Proof
Let $g\in G$ and $a\in\mathcal A$ be a positive operator. Then by assumption (A) and Remark \ref{abscont}
$$0\leq \varphi(g(a))=\varphi(x_ga)\leq ||x_g||\varphi(a)\leq \lambda \varphi(a).$$
Moreover, one has
\begin{eqnarray*}
\varphi(a)&=&\varphi(g^{-1}(g(a))\\
          &=&\varphi(x_{g^{-1}}(g(a))\\
          &\leq&||x_{g^{-1}}||\varphi(g(a))
\end{eqnarray*}
which proves that
$$\varphi(g(a))\geq \frac{1}{||x_{g^{-1}}||}\varphi(a)\geq \frac{1}{\lambda}\varphi(a).$$
Since $\varphi_g=\varphi\circ g: a\mapsto\varphi(x_ga)$ is a state, the linear form
$a\mapsto \varphi(a (x_g^{-1})^* )=\varphi_g(x_g^{-1}a(x_g^{-1})^*)$ is positive. Then, for all $a\in \mathcal A$, we have
$$\varphi(a (x_g^{-1})^*)=\overline{\varphi(a^* (x_g^{-1})^*)}=\overline{\varphi((x_g^{-1}a)^*)}=\varphi(x_g^{-1}a)=L_{x_g^{-1}}\varphi(a) $$
thanks to the positivity of the above-mentioned linear form and $\varphi$, and we can use Lemma \ref{SZ} again as $L_{x_g^{-1}}\varphi\geq 0$. Similarly $\varphi=L_{x_g}(L_{x_g^{-1}}\varphi)\geq 0$ gives the last inequality.
\EndProof

Our main result states that this assumption is equivalent to have an invariant normal faithful state for the action. We borrow some ideas of \cite{AS}.

\begin{thm}\label{thm1} {\rm If $\varphi$ is a normal $G$-quasi-invariant state on $\mathcal A$ with cocycle $x_g$'s satisfying \eqref{A}, then there exist a normal faithful $G$-invariant state $\psi$ on $\mathcal A$ and a bounded operator $d\in\mathcal A$ such that  
$$\psi(a)=\varphi(da),\quad a\in \mathcal A.$$
Conversely, if there exist a normal faithful $G$-invariant state $\psi$ on $\mathcal A$ and a bounded invertible operator $d\in\mathcal A$ such that for all $g\in G$ 
$$\psi(a)=\varphi(da),\quad a\in \mathcal A,$$
then $\varphi$ is a normal $G$-quasi-invariant state on $\mathcal A$ with cocycle $x_g$'s satisfying \eqref{A}.
}
\end{thm}
\Proof
On $\mathcal A$, define the linear map
$$\Gamma_g(a):=x_{g^{-1}}g(a), \quad a\in \mathcal A.$$
Since $g$ is normal and $x_{g^{-1}}$ is bounded, it follows that $\Gamma_g$ is $\sigma$-weakly continuous. Moreover, the maps $\Gamma_g$ satisfy the following properties:
\begin{enumerate}
\item[(i)] For all $g,h \in G$ , $\Gamma_g(x_h)=x_{hg^{-1}}$.
This is exactly \eqref{eq:cocycle}. 
\item[(ii)] For all $h,g\in G$, one has $\Gamma_{gh}=\Gamma_g\Gamma_h$. 

If $a\in \mathcal A$, computing with \eqref{eq:cocycle}
\begin{eqnarray*}
\Gamma_{gh}(a)&=&x_{(gh)^{-1}}\;gh(a)\\
              &=&x_{h^{-1}g^{-1}}\;gh(a)\\
							&=&x_{g^{-1}}g(x_{h^{-1}})g(h(a))\\
							&=&x_{g^{-1}}g(\Gamma_h(a))\\
							&=&\Gamma_g(\Gamma_h(a)).
\end{eqnarray*}
\item[(iii)]  For all $g\in G$, $\varphi\circ \Gamma_g=\varphi$.

If $a\in \mathcal A$, easy computations give
$$\varphi(\Gamma_g(a))=\varphi(x_{g^{-1}}g(a))=\varphi(g^{-1}(g(a)))=\varphi(a).$$
\item[(iv)] For all $g\in G$ and $a,b\in\mathcal A$, $\Gamma_g(ab)=\Gamma_g(a)(x_{g^{-1}})^{-1}\Gamma_g(b)$.

\item[(v)] For all $g\in G$ and $a\in \mathcal A$
$$(\Gamma_g(a))^*=(x_{g^{-1}})^{-1}\Gamma_g(a^*)(x_{g^{-1}})^{*}.$$
\end{enumerate}

Let $K$ be the $\sigma$-weak closure of $\mbox{conv}\{x_g\}_{g\in G}$. Since the set $\{x_g\}_{g\in G}$ is uniformly bounded, $K$ is convex and compact with respect to the $\sigma$-weak topology.

Moreover, by (i), $\Gamma_g(K)\subset K$ for all $g\in G$. 

Now define a continuous norm $\rho$ on $\mathcal A$ by $\rho(x)=\varphi(xx^*)^{1/2}$.  Let $a\neq b\in K$, thanks to the second estimate in Lemma \ref{lemm} applied to $\Gamma_g(a-b)(\Gamma_g(a-b))^*\geq 0$ for $g^{-1}$:
$$\rho(\Gamma_g(a)-\Gamma_g(b))^ 2\geq 
\frac 1\lambda \varphi(\Gamma_g(a-b)(\Gamma_g(a-b))^*((x_{g^{-1}})^{-1})^*).$$
Using the algebraic identities (v), (iv) and (iii) 
\begin{eqnarray*}
\lambda \rho(\Gamma_g(a)-\Gamma_g(b))^ 2
  &\geq&\varphi(\Gamma_g(a-b)(x_{g^{-1}})^{-1}\Gamma_g((a-b)^*))\\
	&=&\varphi(\Gamma_g((a-b)(a-b)^*))\\
		&=&\varphi((a-b)(a-b)^*)>0.
\end{eqnarray*}

Therefore $\Gamma_G$ is non-contracting and by the Ryll-Nardzewski fixed point theorem \cite{RN} there exists $d\in K$ such that $\Gamma_g(d)=d$, i.e. $d=x_{g^{-1}} g(d)$ for all $g\in G$. Notice that for all $g\in G$, $\varphi(x_g)=1$. Then for all $x\in K$ $$\varphi(x)=1=\varphi(d).$$ 
It follows that the operator $d$ is non zero.
Define
$$\psi(a):=\varphi(da),\quad a\in \mathcal A.$$
Then for all $g\in G$
\begin{eqnarray*}
\psi(g(a))=\varphi(dg(a))=\varphi\big(g(g^{-1}(d)a)\big)=\varphi\big(x_gg^{-1}(d)a\big)=\varphi(da)=\psi(a)
\end{eqnarray*}
which proves that $\psi$ is $G$-invariant. By Lemma \ref{lemm} for any $y\in K$, we have $\frac 1\lambda \varphi\leq  \varphi.y \leq \lambda \varphi$. In particular $\psi$ is faithful as $\varphi$ is.


Conversely, let $d$ be a bounded invertible element of $\mathcal A$ such that 
$$\psi(a):=\varphi(da), \quad a\in \mathcal A$$
is $G$-invariant. It follows that for all $a\in \mathcal A$, $\varphi(a)=\psi(d^{-1}a)$, and for all $g\in G$
\begin{eqnarray*}
\varphi(g(a))&=&\psi(d^{-1}g(a))\\
            &=&\psi(g(g^{-1}(d^{-1})a)) \\
						&=&\psi(g^{-1}(d^{-1})a)\\
						&=&\varphi(dg^{-1}(d^{-1})a).
\end{eqnarray*}
This proves that $\varphi$ is $G$-quasi-invariant state with cocycle 
\begin{equation}\label{cocycle}
x_g=d g^{-1}(d^{-1}).
\end{equation}
Since $d$ is a bounded invertible operator, then for any $g\in G$
$$g^{-1}(d^{-1}(d^{-1})^*)\leq ||d^{-1}||^2\;\;\mbox{and}\;\; dg^{-1}(d^{-1}(d^{-1})^*)d^*\leq ||d||^2||d^{-1}||^2. $$
It follows that for any $g\in G$
$$||x_g||^2=||x_gx_g^*||\leq ||d||^2||d^{-1}||^2. $$
Thus the family $\{x_g\}_{g\in G}$ satisfies \eqref{A}.
\EndProof

Now, let $\varphi$ be a normal faithful $G$-strongly quasi-invariant state on $\mathcal A$ with cocycle $x_g'$s satisfying assumption \eqref{A}. Since $x_g$ is a positive operator, one has
$$x_g\leq \lambda,\quad \forall g\in G. $$
Thanks to \eqref{xg-1}, this is equivalent to
\begin{eqnarray}\label{xg-1>}
x_{g^{-1}}\geq \frac{1}{\lambda},\quad\forall g\in G.
\end{eqnarray}
Therefore it is clear that assumption \eqref{A} is equivalent to 
\begin{equation}\label{equiv}
\frac{1}{\lambda}\leq x_g\leq \lambda, \quad\forall g\in G.
\end{equation}
As a consequence of Theorem \ref{thm1}, we prove the following.
\begin{thm}\label{cor}
{\rm $\varphi$ is a normal faithful $G$-strongly quasi-invariant state with cocycle $x_g$'s satisfying \eqref{A} if and only if there exist a normal faithful $G$-invariant state $\psi$ and a bounded positive invertible operator $d\in\mathcal A$ such that for all $g\in G$ 
$$\psi(a):=\varphi(da),\quad a\in \mathcal A.  $$
Moreover, in this case for all $g\in G$, $[d,\;g(d)]=0$.
}
\end{thm}
\Proof
Since $\varphi$ is a normal faithful $G$--strongly quasi-invariant state with cocycle $x_g$'s satisfying assumption \eqref{A} then by Theorem \ref{thm1}, there exist a normal faithful $G$-invariant state $\psi$ and a bounded operator $d\in\mathcal A$ such that for all $g\in G$ 
$$\psi(a):=\varphi(da),\quad a\in \mathcal A.  $$
Remember that from the proof of Theorem \ref{thm1}, $K$ is the $\sigma$-weak closure of $\mbox{conv}\{x_g\}_{g\in G}$ and the cocycle $x_g$'s satisfy \eqref{equiv}, it follows that $d$ satisfies \eqref{equiv} and $d$ is a positive invertible operator in $\mathcal A$. Notice that in this case $x_g=d g^{-1}(d^{-1}) =x_g^*$ which is equivalent to $[d,\;g(d)]=0$.

Conversely, let $\psi$ be a normal faithful $G$-invariant state 
$$\psi(a):=\varphi(da),\quad a\in \mathcal A$$
for some normal state $\varphi$ on $\mathcal A$, where $d$ is a positive invertible operator in $\mathcal A$ such that for all $g\in G$, $[d,\;g(d)]=0$. Then it is clear that $\varphi$ is a normal faithful $G$--strongly quasi-invariant state on $\mathcal A$ with cocycle $x_g=d g^{-1}(d^{-1}) =x_g^*$. Since $d$ is a positive invertible operator, it is clear that the cocycle $x_g'$s satisfies \eqref{A}. 
\EndProof

\section{Unitary implementation associated to quasi--invariant states }\label{unitary-transformation}

When $\varphi$ is $G$-strongly quasi-invariant state on a von Neumann
algebra $\mathcal A$ (for a given action) with cocycle $x_g$'s, it is proved in \cite{AcDh1} that, in the GNS representation of $\varphi$, there exists a unique unitary representation of $G$  implementing the action by conjugation. We will recover it in this section and show that for $G$-quasi-invariant states, one can always spatially implement the action of $G$ (but we don't know if the map $g\mapsto U_g^*$ is a representation). Note that in \cite{Haag}, it is showed that the group $G$ has always a canonical unitary implementation which is different to the one introduced here. In particular in [Theorem 3.2, \cite{Haag}], with the usual notation from there, the condition $u_gJ=Ju_g$ is not satisfied in our case.

Let $G$ be a group of normal $*$-automorphisms of a von Neumann algebra $\mathcal A$ and $\varphi$ be a normal $G$-quasi-invariant state on $\mathcal A$ with cocycle $x_g'$s. 

We will use Haagerup formalism in this section for the standard form of $\mathcal A$ as in section II of \cite{terp}, this will avoid heavy formulas. We denote by $\mathcal N$ the crossed product $\mathcal A\rtimes_{\sigma^\varphi}\mathbb R$ of $\mathcal A$ by the modular group $\sigma^\varphi$ associated with $\varphi$. We denote by $\theta$ the dual action of $\mathbb R$ on $\mathcal N$. Then $\theta$ satisfies $\theta_s(x)=x\;\;(s\in\mathbb R)$ iff $x\in \mathcal A$. Let $\tau$ be the trace on $\mathcal N$ such that $\tau\circ\theta_s=e^{-s}\tau$ and denote by $L_0(\mathcal N,\;\tau)$ the space of $\tau$-measurable operators $x$ affiliated with $\mathcal N$. Define for $p=1,2$
$$L_p(\mathcal A):=\{x\in L_0(\mathcal N,\;\tau)|\;\;\theta_s(x)=e^{-s/p}x,\;\; s\in\mathbb R \}.$$

There is an identification of $\mathcal A_*$ with $L_1(\mathcal A)$  assigning a positive operator $D_\gamma\in L_1(\mathcal A) $ to a normal linear form $\gamma\in \mathcal A_*$  (Proposition 15 \cite{terp}). The duality is defined through the use of a trace-like functional  $\mbox{tr}: L_1(\mathcal A)\to \mathbb C$ with 
$$\langle D_\gamma,\;a\rangle_{L_1(\mathcal A),\;\mathcal A}=\mbox{tr}( D_\gamma a)=\gamma(a).$$
 Since $\varphi$ is faithful state on $\mathcal A$, its density $D_\varphi\in L_1(\mathcal A)$ satisfies supp$(D_\varphi)=1$ by Proposition 4(2c) in \cite{terp}. Note that, by duality, $G$ acts on $L_1(\mathcal A)$ and we have
 $$\langle g^*D_\varphi,\;a\rangle_{L_1(\mathcal A),\,\mathcal A}=\varphi(g(a)).$$
 
 The space $L_2(\mathcal A)$ is naturally a Hilbert space with scalar product $\langle A,B\rangle=\mbox{tr}{(A^*B)}$. Moreover, by Theorem 36 in \cite{terp}, the GNS representation of $(\mathcal A,\;\varphi)$ is the triple  $(\mathcal A\subset \mathcal N,\;L_2(\mathcal A),\;D_\varphi^{1/2})$. 
 
\begin{lem}\label{4.1}\rm{
There exists a positive invertible operator $a_g$ in $\mathcal A$ such that
\begin{equation}\label{equ}
D_{\varphi\circ g}=D_\varphi^{1/2} a_g^2 D_\varphi^{1/2}=x_g^*D_\varphi.
\end{equation}
}
\end{lem}
\Proof
Since $\varphi$ is normal $G$-quasi-invariant and for all $a\in\mathcal A$ and $g\in G$, \\ $\varphi(x_ga)=\varphi(ax_g^*)$, it follows that
\begin{eqnarray}\label{xgDvarphi}
g^*D_\varphi:=D_{\varphi\circ g}=x_g^*D_\varphi=D_\varphi x_g.
\end{eqnarray}
On the other hand, for any positive operator $a\in\mathcal A$ and any $g\in G$, one has
$$\varphi(g(a))=\mbox{tr}{(D_{\varphi\circ g}\,a)}=\varphi(x_g a)\leq ||x_g||\varphi(a)= ||x_g||\,\mbox{tr}{(D_\varphi^{1/2}aD_\varphi^{1/2})}.$$
By Theorem 7 \cite{terp}, we have $0\leq D_{\varphi\circ g}\leq \|x_g\| D_{\varphi}$. As $D_\varphi$ has full support, there exists a positive element $a_g$ of $\mathcal N$ such that
\begin{eqnarray*}
D_{\varphi\circ\, g}=D_\varphi^{1/2}a_g^2D_\varphi^{1/2}.
\end{eqnarray*}
Considering the dual action yields that $\theta_s(a_g)=a_g$ for all $s\in\mathbb R$ and $a_g$ is a positive operator in $\mathcal A$. Now, notice that  by the proof of Lemma \ref{lemm}, $\varphi\leq ||x_{g^{-1}}||\,\varphi\circ g$.
Hence one gets
$$D_{\varphi\circ g}=D_\varphi^{1/2}a_g^2D_\varphi^{1/2}\geq \alpha D_\varphi,$$
where $\alpha=||x_{g^{-1}}||^{-1}$. This proves that $a_g^2\geq\alpha$ and $a_g$ is invertible. 
\EndProof

Actually, the equation $x_g^*D_\varphi=D_\varphi x_g$ implies that
$x_g$ sits in the domain of $\sigma^\varphi_{-i}$ and $\sigma^\varphi_{-i}(x_g)=x_g^*$ (see Theorem 3.25 in \cite{Ta}). In particular, $x_g$ sits in the domain of  $\sigma^\varphi_{-i/2}$ and one can check that
\begin{equation} \label{ag2}
a_g^2=\sigma^\varphi_{-i/2}(x_g).
\end{equation}

Now we prove the following.
\begin{thm}\label{unit-trans}\rm{
For all $g\in G$, the map
\begin{eqnarray*}
U_g:&&\,\mathcal A\,D_\varphi^{1/2}\subset L_2(\mathcal A)\longrightarrow L_2(\mathcal A)\\
&&xD_\varphi^{1/2}\longmapsto g^{-1}(x)D_\varphi^{1/2}a_g
\end{eqnarray*}
extends to a unitary transformation on $L_2(\mathcal A)$ (still denoted by $U_g$). Moreover for any $x\in \mathcal A$
$$U_g^*xU_g=g(x).$$
}
\end{thm}
\Proof
For all $x,y\in\mathcal A$, one has
\begin{eqnarray*}
\langle U_g(xD_\varphi^{1/2}),\,U_g(yD_\varphi^{1/2})\rangle_{L_2(\mathcal A)}&=&\mbox{tr} {(a_g D_\varphi^{1/2}g^{-1}(x)^*g^{-1}(y)D_\varphi^{1/2}a_g)}\\
&=&\mbox{tr} {(D_\varphi^{1/2}a_g^2 D_\varphi^{1/2}g^{-1}(x)^*g^{-1}(y))}\\
&=&\varphi(x^*y)\\
&=&\langle xD_\varphi^{1/2},\,yD_\varphi^{1/2}\rangle_{L_2(\mathcal A)}.
\end{eqnarray*}
Since $\mathcal A\,D_\varphi^{1/2}$ is dense in $L_2(\mathcal A)$, $U_g$ can be extended to an isometry on $L_2(\mathcal A)$. Now, let $h\in L_2(\mathcal A)\cap \mbox{Ran}(U_g)^\bot$. Then, for all $a\in \mathcal A$ 
$$0=\langle h,\,U_g(aD_\varphi^{1/2}) \rangle_{L_2(\mathcal A)}=\mbox{tr}{(h^*g^{-1}(a)D_\varphi^{1/2}a_g)}=\mbox{tr}{(a_gh^*g^{-1}(a)D_\varphi^{1/2})}.$$ 
Notice that $\mathcal A\,D_\varphi^{1/2}$ is dense in $L_2(\mathcal A)$ and $a_g$ is invertible. Therefore, one gets $a_gh^*=0$ and $h=0$. Hence, $U_g$ is a unitary operator on $L_2(\mathcal A)$. 

Now, let $x,y,z\in\mathcal A$. Then, one has
\begin{eqnarray*}
\langle U_g^*xU_g(yD_\varphi^{1/2}),\,zD_\varphi^{1/2}\rangle_{L_2(\mathcal A)}&=&\langle xg^{-1}(y)D_\varphi^{1/2}a_g,\;g^{-1}(z)D_\varphi^{1/2}a_g\rangle_{L_2(\mathcal A)}\\
&=&\mbox{tr}{(a_gD_\varphi^{1/2}(g^{-1}(y))^*x^*g^{-1}(z)D_\varphi^{1/2}a_g )}\\
&=&\varphi(y^*g(x)^*z)\\
&=&\langle g(x)yD_\varphi^{1/2},\;zD_\varphi^{1/2}\rangle_{L_2(\mathcal A)},
\end{eqnarray*}
which proves that $U_g^*xU_g=g(x). $
\EndProof

In the setting of Theorem \ref{thm1}, we have that there exists an invertible $a\in \mathcal A$ such that $D_\psi=D_\varphi^{1/2}a^2 D_\varphi^{1/2}$. Thus there is a partial isometry $v$ such that $D_\psi^{1/2}=va D_\varphi^{1/2}$. Since both operators have full support, $v$ is unitary and  we deduce that there exists an invertible $\gamma\in \mathcal A$ such that $D_\psi^{1/2}=\gamma D_\varphi^{1/2}$. Thus we get $D_\psi=d^*D_\varphi=\gamma D_\varphi \gamma^*$. As above, we get that $\gamma^*$ sits in the domain of $\sigma^\varphi_{-i}$ and $\sigma^\varphi_{-i}(\gamma^*)=\gamma^{-1} d^*$. We conclude that $d^*=\gamma\sigma^\varphi_{-i}(\gamma^*)$. In the same way, let $\gamma_g= g^{-1}(\gamma^{-1})\gamma$ for $g\in G$, then one can also check that $\gamma_g^*$ sits in the domain of $\sigma^\varphi_{-i}$ and that $x_g^*=\gamma_g\sigma^\varphi_{-i}(\gamma_g^*)$.

\begin{prop} \label{represent}\rm{ If $\varphi$ is $G$-strongly quasi-invariant then $g\mapsto U_g^*$ is a unitary representation of the group $G$.
}
\end{prop}
\Proof
Notice that if $\varphi$ is $G$-strongly quasi-invariant, it is showed in \cite{{AcDh1}} that the cocycle $x_g,\;g\in G$, are in the centralizer of $\varphi$. Therefore by \eqref{ag2} it follows that $a^2_g=x_g=x_g^*$ (i.e. $a_g=x_g^{1/2}$). Moreover from \eqref{xgDvarphi}, one has
$$a_gD_\varphi^{1/2}=D_\varphi^{1/2}a_g$$
and hence for all $x\in\mathcal A$, $g\in G$
$$U_g(xD_\varphi^{1/2})=g^{-1}(x)x_g^{1/2}D_\varphi^{1/2}.$$
By using the cocycle property, for any $g,h\in G$, $x\in\mathcal A$, one has
\begin{eqnarray*}
U_gU_h(xD_\varphi^{1/2})=U_g(h^{-1}(x)x_h^{1/2}D_\varphi^{1/2})=g^{-1}(h^{-1}(x))g^{-1}(x_h^{1/2})x_g^{1/2}D_\varphi^{1/2}=U_{hg}(xD_\varphi^{1/2}).
\end{eqnarray*}
It follows that $g\mapsto U_g^*=U_{g^{-1}}$ is a unitary representation of the group $G$.
\EndProof

\section{A Theorem of Kov\'acs and Sz\"ucs for quasi--invariant states with uniformly bounded cocycle}
By Theorem \ref{thm1}, if $\varphi$ is a normal faithful $G$-quasi--invariant state on a von Neumann algebra $\mathcal A$ with cocycle satisfying assumption \eqref{A}, then there exist a normal faithful $G$-invariant state $\psi$ on $\mathcal A$ and a bounded operator $d\in\mathcal A$ such that for all $g\in G$ 
$$\psi(a):=\varphi(da),\quad a\in \mathcal A.$$
Since $\psi$ is a normal faithful $G$-invariant state on $\mathcal A$, it follows from \cite[Theorem 1]{DKS} that there exists a normal faithful $G$-invariant Umegaki conditional expectation $\Phi:\mathcal A\rightarrow \mathcal B$, where $\mathcal B$ is a sub-von Neumann algebra of $\mathcal A$, such that
$$\psi(.)=(\psi|\mathcal B)\circ \Phi(.)  $$
If moreover $d$ is invertible, one gets
\begin{equation}\label{Kovacs-szucs}
\varphi=(\varphi(d.)|\mathcal B)\circ \Phi(d^{-1}.)  
\end{equation}

Now, if $\varphi$ is a $G$-strongly quasi-invariant state and if $g\mapsto U_g^*$ is the unitary representation of $G$ defined by Theorem \ref{unit-trans}, we will show that the conditional expectation $\Phi$ in \eqref{Kovacs-szucs} can also be characterized by $U_g,\;g\in G$ . Define
\begin{equation}\label{*-auto-ug}
u_{g}(a) := U^*_gaU_g=g(a)
\ ;\quad\forall\,g\in G\ ;\quad\forall\,a\in\mathcal{A}
\end{equation}
Let $E_0$ be the orthogonal projection onto the subspace of all $\xi\in L_2(\mathcal A)$ invariant under all $U_g$. We denote by $\mathcal B$ the set of fixed points of $\mathcal A$ under the action of  the automorphisms $u_g$, $g\in G$. Let $F_0=[\mathcal B'E_0]$ be the orthogonal projection onto $\overline{\mbox{span}(\mathcal B'E_0 L_2(\mathcal A))}$. Finally, let $\mathcal M$ be a Godement mean over $G$. The following result follows from a theorem of  Kov\'acs and Sz\"ucs (see \cite{KS} or \cite{AS} and \cite[Theorem 1]{DKS}).

\begin{thm} {\rm Assume that $\varphi$ is a faithful normal $G$-strongly quasi-invariant state with uniformly bounded cocycle (i.e. cocycle satisfying \eqref{A}). Then
\begin{enumerate}
\item[(1)] There exists a unique normal faithful $G$-invariant conditional expectation\\ $\Phi:\mathcal A\rightarrow \mathcal B$. Moreover, $\Phi(b)$ can be defined (equivalently) as
\begin{enumerate}
\item[(i)] the unique element of $\mathcal A$ such that $\Phi(b)E_0=E_0bE_0$
\item[(ii)] the unique element of $\mathcal A$ such that
$$(\xi_1,\;\Phi(b)\xi_2)=\mathcal M\{(\xi_1,\;U_{g^{-1}}bU_g\xi_2)\}$$
\end{enumerate}
\item[(2)] A faithful normal state $w$ on $\mathcal A$ is $G$-strongly quasi-invariant with uniformly bounded cocycle if and only if there exists a positive invertible bounded operator $d$ such that $w(.)=(w(d.)|\mathcal B)\circ \Phi(d^{-1}.)$
\end{enumerate}
}
\end{thm}

\Proof Let $\varphi$ is a faithful normal $G$-strongly quasi-invariant state with uniformly bounded cocycle. By applying Theorem \ref{cor}, there exists a bounded positive invertible operator $d$ and a state a $G$-invariant state $\psi$ such that
$$\psi(a)=\varphi(da)=\langle d^{1/2}D_\varphi^{1/2}, ad^{1/2}D_\varphi^{1/2}\rangle_{L_2(\mathcal A)} \quad \forall a\in \mathcal A.$$
Remember that $x_{g}=d g^{-1}(d^{-1})=x_g^*$. Moreover for all $g\in G$
$$U_gd^{1/2}D_\varphi^{1/2}=g^{-1}(d^{1/2})x_{g}^{1/2}D_\varphi^{1/2}=d^{1/2}D_\varphi^{1/2}.$$
Let $p$ be the orthogonal projection onto $\overline{\mbox{span}(\mathcal B'd^{1/2}D_\varphi^{1/2}})$. Then $p\in\mathcal B$ and $p\leq F_0$. Moreover, one has
\begin{eqnarray*}
\psi(p)&=&\langle d^{1/2}D_\varphi^{1/2}, pd^{1/2}D_\varphi^{1/2}\rangle_{L_2(\mathcal A)}\\
&=&\langle d^{1/2}D_\varphi^{1/2}, d^{1/2}D_\varphi^{1/2}\rangle_{L_2(\mathcal A)}\\
&=&\varphi(d)=\psi(I)=1.
\end{eqnarray*}
Since $\psi$ is faithful, $p=I$ and $F_0=I$. Finally the results of the above theorem follows from \cite[Theorem 1]{DKS}.
\EndProof

Now we provide a commutative easy  example of faithful normal $G$-quasi-invariant state such that the cocycle is not uniformly bounded. Moreover, we show in this example that there is no $G$-invariant conditional expectation.\smallskip\\
{\bf Example}. On $\mathcal A= L^\infty(\mathbb R)$, define the state
$$\varphi(f)=\frac 1 \pi \int_\mathbb R \frac{f(s)}{1+s^2}ds.$$
 The additive group $\mathbb R$ acts on $\mathcal A$ by translation:
$$\tau_t(f)(s):=f(s-t)$$
It is clear that $\varphi$ is $\mathbb R$-quasi-invariant state with bounded cocycle
$$x_t(s)=\frac{1+s^2}{1+(s+t)^2}.$$
The cocycle $(x_t)_t$ does not satisfy assumption (A)
since $\|x_t\|, \|x_t^{-1}\|\geq 1+t^2$. The fixed point of the action are constant functions. Thus every conditional expectation $\Phi$ on $L^\infty(\mathbb R)$ has to  be of the form 
\begin{eqnarray}\label{CE}
\Phi(f)=\psi(f)1 
\end{eqnarray}
for some normal state $\psi$. But there are no functions in $L_1$ that are $\mathbb R$-invariant. 

This example also shows that it seems unlikely to improve Theorem \ref{thm1}.

\section{Strongly quasi-invariant states and semifinite trace}

In this section we extend the result developped by St{$\o$}rmer \cite{St} to the case of strongly quasi-invariant states with uniformly bounded cocycles. Given a semifinite von Neumann algebra $\mathcal A$, a group $G$ of $*$-automorphisms is said to act ergodically on the center $\mathcal Z=\mathcal A\cap \mathcal A'$ of $\mathcal A$ if $\mathcal F(G)\cap \mathcal Z=\mathbb C 1$, where $\mathcal F(G)$ is the set of the fixed points of $G$ in $\mathcal A$. 

\begin{thm} {\rm
Let $\mathcal A$ be a semifinite von Neumann algebra and let $G$ be an amenable
group of $*$-automorphisms of $\mathcal A$ and assume that $G$ acts ergodically on the center
$\mathcal Z$ of $\mathcal A$. Suppose that $\varphi$ is a faithful normal $G$-quasi-invariant state on $\mathcal A$ with cocycle $x_g$'s  satisfying assumption $(A)$. Then, there exists up to a scalar multiple a unique faithful normal $G$-invariant semifinite trace $\tau$ of $\mathcal A$. Furthermore the density $c\in L_1(\mathcal A,\tau)^+$ of $\varphi$ with respect to $\tau$ satisfies
\begin{equation}\label{g-1c}
(g^{-1})^*(c)=cx_{g^{-1}},\quad x_g^*c=cx_g, \quad \forall g\in G,
\end{equation}
where $(g^{-1})^*:L_1(\mathcal A,\tau)\rightarrow L_1(\mathcal A,\tau)$ is the predual map of $g^{-1}$ on $\mathcal A$.
}
\end{thm}

\Proof
Under assumptions of the above theorem, it follows from Theorem \ref{thm1} that there exists a normal faithful $G$-invariant state $\psi$.
Applying [\cite{St}, Theorem 1] for $\psi$, there exists up to a scalar multiple a unique faithful normal $G$-invariant semifinite
trace $\tau$ on $\mathcal A$.
Let $c\in L_1(\mathcal A,\tau)^+$ be the density of $\varphi$. For all $a\in\mathcal A$, we have
$$\tau((g^{-1})^*(c)a)=\tau(cg^{-1}(a))=\varphi(g^{-1}(a))=\varphi(x_{g^{-1}}a)=\tau(cx_{g^{-1}}a). $$
This proves the first equality \eqref{g-1c}. For the second, note that for all $a\in\mathcal A$
\begin{eqnarray*}
\tau(cx_ga)=\varphi\circ g (a)=\overline{\varphi\circ g (a^*)}=\overline{\tau(cx_ga^*)}=\tau(x_g^*ca).
\end{eqnarray*}
This implies that $x_g^*c=cx_g$ in $L_1(\mathcal A,\tau)$.
\EndProof

\noindent {\bf Example}. Taking again $\mathcal A= L^\infty(\mathbb R)$ and $G$ to be the $ax+b$ group of all automorphisms of  $\mathcal A$ of the form $\pi_{(a,b)}f(t)=f(at+b)$ for $(a,b)\in \mathbb R^*\times \mathbb R$.  Then $G$ is amenable (since it is solvable), it admits a cocycle for the state given by the measure $\frac{1}{\pi}\frac 1{1+t^2}dt$; $x_{(a,b)}(t)=\frac {a(1+t^2)}{a^2+(t-b)^2}$ that is not uniformly bounded. Moreover, one can check that the fixed point algebra by $G$ consists only in constant functions. Similarly, one also easily sees that there are no measures on $\mathbb R$ invariant by $G$. Thus one {cannot} remove the assumption $(A)$ even in a commutative situation.

\section*{Acknowledgement} A. Dhahri is a member of GNAMPA-INdAM and he has been supported by the MUR grant Dipartimento di Eccellenza 2023-2027 of Dipartimento di Matematica, Politecnico di Milano.

\section*{Statements and declarations}
\textbf{Ethical Approval}: Not applicable.\\
\textbf{Competing interests}: The authors declare that they have no conflict of interest.\\
\textbf{Data availability statement}: The manuscript has no associated data.

 \end{document}